\documentclass[11pt]{amsart}
\newtheorem{theorem}{Theorem}
\newtheorem{lemma}{Lemma}
\newtheorem{prop}{Proposition}
\newtheorem{coroll}{Corollary}
\theoremstyle{definition}
\newtheorem{example}{Example}
\newtheorem{rem}{Remark}
\newtheorem{quest}{Question}
\renewcommand{\int}{\operatorname{int}}
\newcommand{\cl}{\operatorname{cl}}
\newcommand{\st}{\operatorname{st}}

\begin{document}
\setlength{\unitlength}{0.01in}
\linethickness{0.01in}
\begin{center}
\begin{picture}(474,66)(0,0)
\multiput(0,66)(1,0){40}{\line(0,-1){24}}
\multiput(43,65)(1,-1){24}{\line(0,-1){40}}
\multiput(1,39)(1,-1){40}{\line(1,0){24}}
\multiput(70,2)(1,1){24}{\line(0,1){40}}
\multiput(72,0)(1,1){24}{\line(1,0){40}}
\multiput(97,66)(1,0){40}{\line(0,-1){40}}
\put(143,66){\makebox(0,0)[tl]{\footnotesize Proceedings of the Ninth Prague Topological Symposium}}
\put(143,50){\makebox(0,0)[tl]{\footnotesize Contributed papers from the symposium held in}}
\put(143,34){\makebox(0,0)[tl]{\footnotesize Prague, Czech Republic, August 19--25, 2001}}
\end{picture}
\end{center}
\vspace{0.25in}
\setcounter{page}{15}
\title[Maximal $G$-compactifications of $G$-spaces]{The maximal 
$G$-compactifications of $G$-spaces with special actions}
\author{V.~A.~Chatyrko}
\address{Department of Mathematics\\
Link\"oping University\\
581 83 Link\"oping\\
Sweden}
\email{vitja@mai.liu.se}
\author{K.~L.~Kozlov}
\address{Department of Mechanics and Mathematics\\
Moscow State University\\
117234 Moscow\\
Russia}
\thanks{The second author is supported by RFFI, project 99-01-00128}
\thanks{V.~A.~Chatyrko and K.~L.~Kozlov,
{\em The maximal $G$-compactifications of $G$-spaces with special actions},
Proceedings of the Ninth Prague Topological Symposium, (Prague, 2001),
pp.~15--21, Topology Atlas, Toronto, 2002}
\begin{abstract}
An action on a $G$-space induces uniformities on the phase space.
It is shown when the maximal $G$-compactification of a $G$-space
can be obtained as a completion of the phase space with respect
to one of these uniformities. Structure of $G$-spaces with
special actions is investigated.
\end{abstract}
\maketitle

This paper is the continuation of the previous work of the 
authors~\cite{CK} and is partially supported by Kungliga Vetenskapademian,
project 12529. 
Besides old results which are now proved using another techniques the new
ones are presented.

All spaces are assumed to be Tychonoff and mappings are continuous
mappings of spaces.
Let $R$ denote the real numbers, and nbd is an abridged notation for
neighbourhood.

Let $G$ be a topological group.
By a $G$-space $X$ we mean a Tychonoff space $X$ (phase space) with a
continuous action of group $G$.
If for a $G$-space $X$ there exist a compact $G$-space $bX$ and an 
equivariant dense embedding of $X$ into $bX$ then we call $bX$ a 
$G$-compactification (see, for example, \cite{V}).
If a $G$-space has a $G$-compactification then (see, for example,
\cite{AS}) there is the largest element $\beta_G X$ among all
$G$-compactifications which is called the maximal $G$-compactification.

Uniform structures are introduced by coverings~\cite{E}, and we say that
the uniformity $U_1$ is finer than the uniformity $U_2$ if 
$U_2\subset U_1$.
If the topology of a topological space and the topology induced by a
uniformity on it are the same then we say that the uniformity is
compatible with the topology of the space.
If $X$ is a $G$-space then the uniformity $U$ on $X$ is called {\it
invariant} if for any $g\in G$ and $\gamma\in U$ $g\gamma\in U$.

In 1975 J.\,de~Vries~\cite{V} introduced the notion of a {\it bounded 
action} 
(an action on a space $X$ is bounded if there exists a uniformity $U$
on $X$ compatible with its topology such that for any $u\in U$ there is a
nbd $O$ of identity in $G$ such that the pair of points $x$ and $g x$
belong to one element of $u$ for any $x\in X$ and any $g\in O$) 
and proved that a $G$-space $X$ has a $G$-compactification iff the action
is bounded.

In 1984 M.\,G.~Megrelishvili~\cite{M} introduced the concept of an {\it
equiuniformity} 
(uniformity is equiuniformity if it is compatible with the topology of the
phase space, invariant and the action is bounded by it) 
and proved the following theorem.

\newtheorem*{thmA}{Theorem A}

\begin{thmA}
(M.\,G.~Megrelishvili~\cite{M})
If $U$ is an equiuniformity on a $G$ space $X$ then its completion
$\tilde X$ with respect to $U$ is a $G$ space.
Besides if $f:X\to Y$ is an equivariant uniformly continuous mapping to a
complete uniform space $Y$ then there exists the unique equivariant
continuous mapping $\tilde f:\tilde X\to Y$ such that $\tilde f\circ i=f$
where $i$ is a natural embedding of $X$ into $\tilde X$.
\end{thmA}

The following statements are evident.

\begin{prop}
Let $U_\alpha, \alpha\in A,$ be the family of equiuniformities on a $G$ 
space $X$. 
Then its least upper bound is an equiuniformity.
\end{prop}

\begin{coroll}
Among all equiuniformities there is a maximal one.
\end{coroll}

\begin{prop}
If $U$ is an equiuniformity on a $G$ space $X$ then the set of all
coverings which can be refined by a finite covering from $U$ is an
equiuniformity.
\end{prop}

\begin{coroll}
Among all equiuniformities there is a maximal totally bounded one.
\end{coroll}

Let $A$ be the family of all open nbds of the identity in $G$.
Every $O\in A$ sets two coverings of a $G$-space $X$:
$$\gamma_O=\{ Ox: x\in X\}\ \mbox{and}\ \bar\gamma_O=\{\cl (Ox): x\in X\}.$$
Denote by $U_G$ ( $\bar U_G$ ) the family of all coverings of $X$ which
have a refinement of the form $\gamma_O$ ($\bar\gamma_O$), $O\in A$.
It may be easily checked that the family $U_G$ is a uniformity on $X$ (not
nessesary compatible with the topology of $X$).

\begin{rem}
If $\bar U_G$ is a uniformity on $X$ then the uniformity $U_G$ is finer
than $\bar U_G$, but they may not be compatible with the topology of the
phase space.
\end{rem}

Now we shall reformulate J.\,de~Vries's criterion mentioned above.

\newtheorem*{thmB}{Theorem B}

\begin{thmB}
(J.\,de~Vries~\cite {V})
A $G$-space $X$ has a $G$-compactification iff there is a uniformity $U$
on $X$ compatible with its topology such that $U_G$ is finer than $U$.
\end{thmB}

Let $U^*$ be the totally bounded uniformity on the space $X$ compatible 
with its topology such that any bounded continuous function on $X$ is 
uniformly continuous with respect to it.
It is the maximal totally bounded uniformity on $X$.

\begin{theorem} 
Let $X$ be a $G$-space.
If the uniformity $U_G$ is finer than $U^*$ then
$$\beta_G X=\beta X.$$
\end{theorem}

\begin{proof}
It is easy to see that $U^*$ is an equiuniformity and
the rest follows from Corollary 2 and Theorem A.
\end{proof}

\begin{lemma}
Let $U$ be a uniformity on a $G$-space $X$ compatible with its topology.
If $\bar U_G$ is a uniformity on $X$ then the following conditions are
equivalent:
\begin{itemize}
\item[{\rm (1)}] $U_G$ is finer than $U$;
\item[{\rm (2)}] $\bar U_G$ is finer than $U$.
\end{itemize}

Moreover, if the uniformity $U_G$ is compatible with the topology of $X$
then uniformities $U_G$ and $\bar U_G$ are the same.
\end{lemma}

\begin{proof}
Let the uniformity $U$ be generated by the family of coverings.
In order to show $(1)\Rightarrow (2)$ for any $v\in U$ take $v'\in U$ such
that $v'$ is a star refinement of $v$. Then the covering $[v']\in U$ which
consists of closures of elements of $v'$ is a refinement of $v$. 
Take $\gamma_O\in U_G$ such that $\gamma_O$ is a refinement of $v'$.
Then $\bar\gamma_O\in\bar U_G$ is a refinement of $[v']$.
From this it follows that $\bar\gamma_O$ is a refinement of $v$ and hence
$\bar U_G$ is finer than $U$.

The implication $(2)\Rightarrow (1)$ follows from Remark 1.

From Remark 1 it follows that $U_G$ is finer than $\bar U_G$. 
If $U_G$ is compatible with the topology of $X$ then instead of $U$ we can
take $U_G$ in our lemma. 
Then $\bar U_G$ is finer $U_G$ also. 
Hence $U_G$ and $\bar U_G$ are the same.
\end{proof}

\begin{prop}
If $\bar U_G$ is a uniformity compatible with the topology of $X$ then
it is a maximal equiuniformity.
\end{prop}

\begin{proof}
Since the uniformity $U_G$ is finer than $\bar U_G$, it follows from
Theorem B that the action is bounded by the
uniformity $\bar U_G$.

In order to prove that the uniformity $\bar U_G$ is
invariant it is sufficient to show that for any
$\bar\gamma_O=\{\cl(Ox): x\in X\}$, $O\in A$, $g\bar\gamma_O\in\bar U_G$.
Since for any $g\in G$ the mapping $g: X\to X, x\to gx$ is
a homeomorphism it follows that $g(\cl(Ox))=\cl((gO)x)$.
Take $U=gOg^{-1}$. Then $U\in A$ and $(gO)x=(Ug)x$ for any
$x\in X$. Thus $g\bar\gamma_O=\bar\gamma_U$.

Hence $\bar U_G$ is an equiuniformity. Its maximality
follows from Lemma 1.
\end{proof}

The proof of the following theorem immediately follows from Proposition 3,
Corollary 2 and Theorem A.

\begin{theorem} 
Let $X$ be a $G$-space.
If $\bar U_G$ is a uniformity compatible with the topology of $X$ then
$$\beta_G X\ \mbox{is the Samuel compactification of}\ X\ \mbox{with
respect to}\ \bar U_G.$$
\end{theorem}

The next example shows that the usage of uniformity $\bar U_G$ gives us
more opportunities in finding maximal $G$-compactifications.

\begin{example}
Let $S=\{ z\in C: |z|=1\}$ be a unit circle on the complex plain, and $a$
be such an element of $S$ that $a^n\ne 1, n\in N$.
Let us put $G=\{a^n: n\in Z\}$ (it is a group with a natural
multiplication), $X=S\setminus G$ and the action of $G$ on $X$ is induced
by multiplication in $C$.

The uniformity $U_G$ is not compatible with the topology of $X$
because the group $G$ is countable and the cardinality of each
nonempty open set of $X$ is uncountable
and the uniformity $\bar U_G$ is compatible
because $G$ is a dense subset of $S$.
\end{example}

\begin{rem}
Earlier Theorems {\rm 1} and {\rm 2} were proved in another way 
in~{\rm\cite{CK}} using results of J.~de Vries~{\rm\cite{V}} and
Yu.\,M.~Smirnov~{\rm\cite{AS}}.
\end{rem}

We can characterize the case when $\bar U_G$ is the uniformity compatible
with the topology of the phase space.

\begin{theorem} 
Let $X$ be a $G$-space.
The family $\bar U_G$ is a uniformity compatible with the topology of $X$
iff the action has the property:
\begin{itemize}
\item[{\rm (a)}] 
for any $x\in X$ and any nbd $O\in A$ there exists $y\in X$ such that 
$x\in\int \cl (Oy)$.
\end{itemize}
\end{theorem}

\begin{proof}
First of all let us notice that property (a) is equivalent
to the following one:
$$\mbox{the family}\ \{\int \cl (Ox): x\in X\}\ \mbox{is a covering of}\
X\ \mbox{for any nbd}\ O\in A.$$

Since the universal uniformity is finer than $\bar U_G$ an open covering
of $X$ may be refined in any covering $\{ \cl (Ox): x\in X\}$ from 
$\bar U_G$. 
From this necessity immediately follows.

In order to prove sufficiency we must first of all check that $\bar U_G$
is the uniformity (see, for example, \cite[page 524]{E}).
Recall that $A$ is the family of all open nbds of identity in $G$.

\begin{itemize}
\item[1.] Right from the definition of $\bar U_G$ it follows that if
$\gamma\in\bar U_G$ and $\gamma$ is a refinement of a covering
$\beta$ of $X$ then $\beta\in\bar U_G$.

\item[2.] It is evident that if $\beta_1$ and $\beta_2\in\bar U_G$
be such that $\bar\gamma_V$ and $\bar\gamma_W$ are refined in
$\beta_1$ and $\beta_2$ for some $V, W\in A$ respectively, then
for $O\in A$ such that
$O\subset V\cap W$ we have that $\bar\gamma_O$ is refined both in
$\bar\gamma_V$ and $\bar\gamma_W$ and hence in $\beta_1$ and $\beta_2$.

\item[3.] For $\beta\in\bar U_G$ let $V\in A$ be
such that $\bar\gamma_V$ is refined in $\beta$. Take $O\in A$
such that $O=O^{-1}$ and $O^3\subset V$. We shall prove
that $\bar\gamma_O$ is a barycentric refinement of $\bar\gamma_V$.

Let us show that $O\cl (W x)\subset\cl (O W x)$ for any nbds $O$ and $W$
of identity in $G$. If $a\in O\cl (W x)$ then
$a=h t$, where $h\in O$ and $t\in\cl (W x)$. Since the action is
continuous for any nbd $V_a$ of $a$ there are a nbd $V_t$ of $t$
such that $h V_t\subset V_a$. Thus there is
$t'\in V_t\cap W x$ such that $ht'\in V_a$.
Hence, $a\in\cl (O W x)$.

For any $x\in X$ there exists $z\in X$ such that $x\in\int \cl (Oz)$.
Now if $x\in\cl (Oy)$ then $\cl (Oz)\cap O y\ne\emptyset$ since
$\int \cl (O z)$ is a nbd of $x$. From this it follows that
$$y\in O^{-1}\cl (O z)\subset\cl (O^2 z)\ \mbox{and}\
O y\subset O\cl (O^2 )\subset\cl (O^3 z)\subset\cl (V z).$$

Thus $\cl (O y)\subset\cl (V z)$ and so $\st (x,\bar\gamma_O)\subset
\cl (V_z)$. Hence, $\bar\gamma_O$ is a barycentric refinement of
$\bar\gamma_V$.

Using the same process, we can find a barycentric refinement of
$\bar\gamma_O$ which would be the star refinement of
$\bar\gamma_V$~\cite[Lemma 5.1.15]{E}.

\item[4.] Let $x, y$ be a pair of distinct points of $X$.
Since the action is continuous and $X$ is a Tychonoff space
there are nbd $O\in A, O^{-1}=O$ and nbds $W_x, W_y$ of $x$ and
$y$ respectively, such that $\cl (O W_x)\cap\cl (O W_y)=\emptyset$.
Let us show that no element of the cover $\bar\gamma_O$ contains
both $x$ and $y$. Indeed, if $x\in\cl (O z)$ and $y\in\cl (O z)$
for some $z\in X$ then $W_x\cap O z\ne\emptyset$ and
$W_y\cap O z\ne\emptyset$. From this it follows that $z\in O^{-1}W_x$
and $z\in O^{-1}W_y$ and, hence, $O W_x\cap O W_y\ne\emptyset$.
This is a contradiction with the choice of nbds $O, W_x$ and $W_y$.
\end{itemize}

So all conditions for the uniformity $\bar U_G$ are fulfilled.

Since for any $x\in X$ and any nbd $O$ of identity in $G$ there exists
$z\in X$ such that $x\in\int \cl (Oz)$
then an open covering can be refined in
any covering from $\bar U_G$. So every open set in topology
induced by uniformity $\bar U_G$ is open in $X$.
If $W$ is open in $X$ and $x\in W$ then there exist $O\in A$ and
a nbd $V$ of $x$ such that $O=O^{-1}$ and $\cl (O^2V)\subset W$.
If $x\in\cl (Oy)$ and $z\in\cl (Oy)$ then there exist
$x_1\in V$ and $h\in O$ such that $x_1=hy$. From this it follows that
$y\in Ox_1$ and $z\in\cl (O^2x_1)\subset\cl (O^2V)\subset W$.
Hence, $\st (x,\bar\gamma_O)\subset W$ and so
$W$ is open in the topology induced be the uniformity.
\end{proof}

\begin{prop}
Consider the following properties for a $G$-space $X$.
\begin{itemize}
\item[{\rm (a)}] 
for any $x\in X$ and any nbd $O\in A$ there exists $y\in X$ such that 
$x\in\int \cl (Oy)$.
\item[{\rm (b)}] 
for any $x\in X$ and any nbd $O\in A$ $x\in\int \cl (Ox)$,
\item[{\rm (c)}] 
for any $x\in X$ and any nbd $O\in A$ $x\in\int (Ox)$,
\end{itemize}
Then {\rm (c)} $\Longrightarrow$ {\rm (b)} $\Longrightarrow$ {\rm (a)}
and the inverse implications are not valid.
\end{prop}

\begin{proof}
The implications (c) $\Longrightarrow$ (b) $\Longrightarrow$ (a) are
evident.

(c) $\not\Longleftarrow$ (b). 
Consider the following example.
Let $Q$ be the set of rational numbers of the interval $I=(0, 1)$.
There is a natural linear order on $Q$.
Let $G$ be a group of all order preserving homeomorphisms of
$Q$~\cite[page 18]{E} with the topology of uniform
convergence~\cite[page 329]{E}. Using Theorem 2 one may show that
$\beta_G Q=[0, 1]$ and
$X=I$ is an invariant subset of $\beta_G Q$. Now it is easy to
see that a $G$-space $X$ satisfies (b) but not (c).

(b) $\not\Longleftarrow$ (a). Consider the following example.
Put $X=\beta_G Q$ where $Q$ and $G$ as above. It is easy to
see that a $G$-space $X$ satisfies (a) but not (b) because the action has
fixed points.
\end{proof}

\begin{rem}
$G$-spaces with property {\rm (b)} were examined by V.\,V.~Uspenski\v{\i}
in~{\rm\cite{U}}.
\end{rem}

Below we shall describe $G$-spaces with properties listed above.

\begin{lemma}
If a $G$-space $X$ satisfies property {\rm (a)} then for any points 
$x, y\in X$ we have either 
$$\int \cl (Gx)=\int \cl (Gy)
\mbox{ or }
\int \cl (Gx)\cap\int \cl (Gy)=\emptyset.$$
\end{lemma}

\begin{proof}
For the proof it is sufficient to show that 
$$\mbox{if }
\int \cl (Gx)\cap\int \cl (Gy)\ne\emptyset
\mbox{ then }
\int \cl (Gx)\subset\cl (Gy).$$

Let $z\in\int \cl (Gx)$ and $O_z$ is an arbitrary nbd of $z$.
We may take $z'\in O_z\cap Gx$, $z''\in\int \cl (Gy)\cap Gx$ and 
$g\in G$ such that $gz''=z'$. 
From the continuity of the action it follows that there exists a nbd
$O_{z''}$ such that $gO_{z''}\subset O_z$.
If we take $y'\in O_{z''}\cap Gy$ then $gy'\in O_z$ and hence
$Gy\cap O_z\ne\emptyset$. Since $O_z$ is an arbitrary nbd of $z$ it
follows that $z\in\cl (Gy)$ and so $\int \cl (Gx)\subset\cl (Gy)$.
\end{proof}

\begin{coroll}
If a $G$-space $X$ satisfies property {\rm (a)} then either 
$$\int \cl (Gx)=\cl (Gx) 
\mbox{ or }
\int \cl (Gx)=\emptyset$$ 
for any point $x\in X$.
\end{coroll}

\begin{theorem}
\mbox{}
\begin{itemize}
\item[{\rm A)}] 
If a $G$-space $X$ satisfies property {\rm (a)} then $X$ is a disjoint
union of clopen sets and each clopen set from this union is the closure
of an orbit of some point and thus it contains the continuous one-to-one
image of some quotient space of group $G$.
\item[{\rm B)}]
If a $G$-space $X$ satisfies property {\rm (b)} then in addition to {\rm
A)} each such clopen set is the closure of orbit of its any point. 
But orbits of different points from the common clopen set may not be 
homeomorphic.
\item[{\rm C)}]
If a $G$-space $X$ satisfies property {\rm (c)} then each such clopen set
is homeomorphic to the orbit of its any point.
\end{itemize}
\end{theorem}

\begin{proof}
Proof of the statements A) and B) follows from Lemma 2 and Corollary 3.
The second example from Proposition 4 shows that in case of action
with property (a) not orbit of any point may be taken
(there may be fixed points). In case of
action with property (b) the first example from Proposition 4 shows
that different orbits (rationals and irrationals) may not be
homeomorphic.

In case of action with property (c) we have for any $x\in X$ an open
mapping $g\to gx$ of $G$ into $X$ and hence a homeomorphism of some
quotient space of $G$ onto its orbit.
\end{proof}

The following questions are not yet known to authors.

\begin{quest} 
When does the family $\bar U_G$ generate uniformity on $X$ (not nessesary
compatible with the topology of $X$)?
\end{quest}

\begin{quest}
Let a $G$-space $X$ satisfy property {\rm (a)}. 
Does there exist a dense invariant subspace $X'$ of $X$ such that the
restriction of action on it has property {\rm (b)}?
\end{quest}

\begin{quest}
Let a $G$-space $X$ satisfy property {\rm (b)}. 
Does there exist a dense invariant subspace $X'$ of $X$ such that the
restriction of action on it has property {\rm (c)}?
\end{quest}

\begin{quest}
Can every compactification of a Tychonoff space be obtained as a
$G$-compactification for some acting group $G$ on $X$?
\end{quest}


\begin{thebibliography}{1}

\bibitem{AS}
S.~A. Antonjan and Ju.~M. Smirnov, \emph{Universal objects and bicompact
  extensions for topological groups of transformations}, Dokl. Akad. Nauk SSSR
  \textbf{257} (1981), no.~3, 521--526. \MR{82i:54072}

\bibitem{CK}
V.~A. Chatyrko and K.~L. Kozlov, \emph{Dimension of maximal equivariant compact
  extensions}, Preprint LiTH-MAT-R-2001-11, Link{\"o}ping University, 2001.

\bibitem{V}
Jan de~Vries, \emph{On the existence of ${G}$-compactifications}, Bull. Acad.
  Polon. Sci. S\'er. Sci. Math. Astronom. Phys. \textbf{26} (1978), no.~3,
  275--280. \MR{58 \#31002}

\bibitem{E}
Ryszard Engelking, \emph{General topology}, PWN---Polish Scientific Publishers,
  Warsaw, 1977, Translated from the Polish by the author, Monografie
  Matematyczne, Tom 60. [Mathematical Monographs, Vol. 60]. \MR{58 \#18316b}

\bibitem{M}
M.~G. Megrelishvili, \emph{Equivariant completions and compact extensions},
  Soobshch. Akad. Nauk Gruzin. SSR \textbf{115} (1984), no.~1, 21--24.
  \MR{86m:54054}

\bibitem{U}
V.~V. Uspenski{\u\i}, \emph{Topological groups and {D}ugundji compact spaces},
  Mat. Sb. \textbf{180} (1989), no.~8, 1092--1118, 1151, translation in Math.
  USSR-Sb. 67 (1990), no. 2, 555--580. \MR{91a:54064}

\end{thebibliography}
\providecommand{\bysame}{\leavevmode\hbox to3em{\hrulefill}\thinspace}
\providecommand{\MR}{\relax\ifhmode\unskip\space\fi MR }
\providecommand{\MRhref}[2]{%
  \href{http://www.ams.org/mathscinet-getitem?mr=#1}{#2}
}
\providecommand{\href}[2]{#2}

\end{document}